\documentclass[10pt]{article}

\usepackage{amsmath}
\usepackage{amssymb}
\usepackage{amsthm,verbatim}
\usepackage{amsfonts}
\usepackage[all]{xy}
\usepackage{mathpazo}
\usepackage{euler}
\usepackage{hyperref}

\DeclareMathOperator{\Spec}{Spec}
\DeclareMathOperator{\ctp}{ctp}
\DeclareMathOperator{\tp}{tp}

\DeclareMathOperator{\Frac}{Frac}
\DeclareMathOperator{\Aut}{Aut}
\DeclareMathOperator{\Gal}{Gal}

\DeclareMathOperator{\WC}{WC}

\newtheorem{theorem}{Theorem}[section]
\newtheorem{example}{Example}[section]

\newtheorem{proposition}{Proposition}[section]
\newtheorem{lemma}{Lemma}[section]
\newtheorem{corollary}{Corollary}[section]
\newtheorem{definition}{Definition}[section]
\newtheorem{remark}{Remark}[section]

\title{Equivariant Zariski Structures}
\author{Vinesh Solanki}

\begin{document}

\maketitle

\begin{abstract}
A new class of noncommutative $k$-algebras (for $k$ an algebraically closed field) is defined and shown to contain some important examples of quantum groups. To each such algebra, a first order theory is assigned describing models of a suitable corresponding geometric space. Model-theoretic results for these geometric structures are established (uncountable categoricity, quantifier elimination to the level of existential formulas) and that an appropriate dimension theory exists, making them Zariski structures.
\end{abstract}

\section{Introduction}
The present paper fits into Zilber's program for constructing novel geometric structures which may be viewed as encoding, in some suitable sense, the geometry of certain noncommutative algebras. What is being done here is a kind of noncommutative algebraic geometry, but from a model-theoretic viewpoint. 
\\
\\
Prior work in this area of model theory by Zilber had concentrated on associating geometric structures to a large class of algebras described as `quantum algebras at roots of unity' (see \cite{Zil06}). The algebras considered there were affine (i.e. finitely generated) $k$-algebras which are large (technically Azumaya) over their centers. The existence of large centers makes such algebras amenable to the techniques of modern algebraic geometry; indeed the corresponding coherent sheaf of algebras over the spectrum of the center of such an algebra functions as a suitable geometric object. The geometric structures associated to quantum algebras at roots of unity in \cite{Zil06} arose quite naturally from such sheaves. 
\\
\\
The equivariant algebras defined in the present paper form a new class of algebras containing significantly more noncommutative objects, e.g. $U_{q}(\mathfrak{sl}_{2}(k))$ where $q$ is not a root of unity. In particular, some of these algebras have very small centers and hence associating the corresponding quasi-coherent sheaf to such an algebra is not geometrically informative. It is demonstrated in this paper that for each equivariant algebra it is possible to associate to it a geometric model-theoretic structure which, in certain favourable cases, has a dimension theory resembling that found for varieties. 
\\
\\
Given that the structures associated to these algebras are of entirely model-theoretic origin, the question is raised as to what criteria one can use to conclude that these structures are sufficiently algebro-geometric. It can be shown that many of the structures of interest considered so far are not interpretable in an algebraically closed field (see \cite{Zil06}, \cite{SSZ14}), hence are not reducible to varieties. Thus for the assigned structures to be suitably algebro-geometric one requires an abstractly given characterization of the geometry of algebraic varieties, but one suitably loose to apply to structures corresponding to various non-commutative algebras. The notion of a Zariski structure (as presented in \cite{Zil10}) fits this purpose aptly:

\begin{definition}
\label{defn:Zar}
Let $\mathbf{X}$ be an infinite set. A \textbf{Zariski structure} on $\mathbf{X}$ consists of a topology on $\mathbf{X}^{n}$ for every $n$ together with an $\mathbb{N}$-valued dimension function $\dim$ on definable subsets of $\mathbf{X}^{n}$ such that the following two collections of axioms are satisfied:
\\
\\
\textbf{Topological axioms}: 
\begin{enumerate}
\item The topology on each $\mathbf{X}^{n}$ is Noetherian.
\item Singletons in $\mathbf{X}^{n}$ are closed. 
\item Cartesian products of closed sets are closed. 
\item The diagonals $x_{i} = x_{j}$ of $\mathbf{X}^{n}$ are closed. 
\item For a tuple $\overline{a} \in \mathbf{X}^{m}$ and any closed set $\mathbf{C} \subseteq \mathbf{X}^{m+n}$ the fiber over $\overline{a}$
\[ \mathbf{C}(\overline{a}, \mathbf{X}^{n}) = \{\overline{b} \in \mathbf{X}^{n}: (\overline{a},\overline{b}) \in \mathbf{C}\} \] is closed. 
\end{enumerate} \textbf{Dimension axioms}: 

\begin{enumerate}
\item The dimension of a point is $0$. 
\item $\dim(\mathbf{S}_{1} \cup \mathbf{S}_{2}) = \max\{\dim \mathbf{S}_{1}, \dim \mathbf{S}_{2}\}$ for all definable subsets $\mathbf{S}_{1}, \mathbf{S}_{2}$. 
\item For $\mathbf{C}$ closed and irreducible in $\mathbf{X}^{n}$ and $\mathbf{C}_{1}$ a closed subset of $\mathbf{C}$, if $\mathbf{C}_{1} \neq \mathbf{C}$ then $\dim \mathbf{C}_{1} < \dim \mathbf{C}$. 
\item For $\mathbf{C}$ irreducible and closed in $\mathbf{X}^{n}$, if $\pi: \mathbf{X}^{n} \rightarrow \mathbf{X}^{m}$ is a coordinate projection map (so $m <n$) then
\[ \dim \mathbf{C} = \dim \pi(\mathbf{C}) + \min_{a \in \pi(\mathbf{C})} \dim (\pi^{-1}(a) \cap \mathbf{C}) \] and there is an subset $\mathbf{V}$ open in $\pi(\mathbf{C})$ such that
\[ \min_{a \in \pi(\mathbf{C})} \dim(\pi^{-1}(a) \cap \mathbf{C}) = \dim(\pi^{-1}(v) \cap \mathbf{C}) \] for every $v \in \mathbf{V}$.
\end{enumerate} 
\end{definition} The dimension function forming part of the data of a Zariski structure can satisfy additional properties. For our purposes, the following notion of {\em presmoothness} is important: it places a bound on how much the dimension of an irreducible component contained in the intersection of two closed irreducible sets can decrease and also gives an abstract characterization of smoothness for algebraic varieties. 

\begin{definition}
\label{defn:presm}
A Zariski structure is said to be \textbf{presmooth} if for any closed irreducible subsets $\mathbf{C}_{1}, \mathbf{C}_{2}$ of $\mathbf{X}^{n}$ the dimension of any irreducible component of $\mathbf{C}_{1} \cap \mathbf{C}_{2}$ is greater than or equal to
\[ \dim \mathbf{C}_{1} + \dim \mathbf{C}_{2} - \dim \mathbf{X}^{n} \] 
\end{definition}

\begin{theorem}
\label{thm:vars:zargeom} An algebraic variety $V$ over an algebraically closed field $k$ in the language containing an $n$-ary relation for each closed subset of $k^{n}$ with $\dim$ given by Krull dimension is a Zariski structure. It is presmooth if $V$ is smooth. 
\end{theorem}

\begin{proof} See \cite{Zil10}, Theorem 3.4.1.
\end{proof} Historically, Zariski structures first appeared (when $\mathbf{X}$ is one-dimensional in some suitable model-theoretic sense) in the paper \cite{HZ96} in a bid to find a class of structures in which the trichotomy conjecture held. Investigations into possible links between Zariski structures and noncommutative algebraic geometry began after the result (also from \cite{HZ96}) that there exist one-dimensional presmooth Zariski structures which cannot be interpreted in an algebraically closed field. Rather than being mathematical pathologies, such structures turned out to correspond naturally to certain noncommutative algebras. In this regard, we mention the paper \cite{SSZ14} as providing an example of such a one-dimensional Zariski structure corresponding to the first Weyl algebra. The techniques developed in \cite{SSZ14} are shown in this paper to be applicable to a larger class of noncommuative algebras.
\\
\\
We now summarize the contents of this paper. We introduce the notion of an equivariant algebra and some examples of interest in Section $2$. In Section $3$, given an equivariant algebra $A$, we associate a structure to it, and this structure is shown to have a first-order axiomatizable theory $T_{A}$. The choice of terminology and the structures considered are motivated by the $G$-equivariant line bundles of geometric representation theory, namely those line bundles $L$ over a variety $V$ endowed with an action of a linear algebraic group $G$, such that 
\[ \mbox{for all $x \in V$ and for all $g \in G$, $g(L_{x}) = L_{gx}$ and $g: L_{x} \rightarrow L_{gx}$ is a linear isomorphism}\] where $L_{x}$ denotes the fiber of $L$ at $x \in V$. The structure associated to an equivariant algebra $A$ is similar: there is an associated (abstract abelian) group $G$ and a variety $V$ corresponding to a commutative subalgebra of $A$; $V$ is endowed with an algebraic action of $G$ and there is an associated bundle of one-dimensional vector spaces over $V$ with the action of $G$ on $V$ inducing linear isomorphisms between these vector spaces. 

Sections $4$ and $5$ and devoted to the model theory of $T_{A}$. Section $4$ contains a complete algebraic characterization of those equivariant $A$ for which models of $T_{A}$ are interpretable in algebraically closed fields. In particular, it is shown that no model of the theory associated to $U_{q}(\mathfrak{sl}_{2}(k))$ is so interpretable. We derive an algebraic characterization of relative category in Section $5$, and for equivariant algebras satisfying this criterion quantifier elimination results are established thus leading to the expected consequences for the category of definable subsets; namely that every definable subset is constructible for an appropriate topology on models. With this topology, an appropriate dimension theory turns each model into a Zariski structure. 
\\
\\
It is worth remarking that at present, it is not clear what use (if any) could be made of the Zariski structure on a given model of $T_{A}$ in order to discern finer structural properties of the equivariant algebra $A$ itself. Nevertheless, stipulating that the associated structure has good model-theoretic properties does place conditions on $A$. It is hoped that such techniques from geometric model theory can function as a prescriptive tool in finding other nice classes of noncommutative algebras. 

\section{Equivariant algebras}
\label{sect:equivalgs}
In this section we define the notion of an equivariant algebra over a field $k$. Firstly, we recall the definition of a {\em skew group ring}.

\begin{definition}
\label{defn:skewgrp} 
Let $R$ be a commutative ring, $G$ a group and suppose that we have a group action 
\[ \varphi_{G}: G \rightarrow \Aut(R) \qquad \varphi_{G}(g): r \mapsto {}^{g} r \] The \textbf{skew group ring} $R * G$ is defined to be the free left $R$-module on generators $\{e_{g}: g \in G\}$ with multiplication defined by
\[ (re_{g})(r'e_{g'}) = r{}^{g}r'e_{gg'} \] 
\end{definition} A skew group ring will in general be non-commutative (it is only commutative if $G$ is abelian and the action of $G$ on $R$ is trivial). Note that there is a natural inclusion of $R$ into $R*G$ by the map $r \mapsto re_{1_{G}}$ where $1_{G}$ is the identity element of $G$ and when dealing with skew group rings we will typically write $r$ instead of $re_{1_{G}}$. 
\\
\\
Now suppose that $R*G$ is a skew group ring for $G$ a finitely generated abelian group and $R$ a commutative $k$-algebra for $k$ a field, with $G$ acting on $R$ by $k$-algebra automorphisms. Then the only source of non-commutativity in $R*G$ comes from the action of $G$ on $R$; this non-commutativity is expressed in the relations $e_{g}r = {}^{g}re_{g}$ for $r \in R$, $g \in G$. In particular, for all $g,g' \in G$, $e_{g}$ and $e_{g'}$ commute. The idea behind an equivariant algebra is to weaken this commutativity restriction on the $e_{g}$ to allow for some relations of the form $e_{g}e_{g'} - e_{g'}e_{g} = r$ (where $r \in R$) but to do so in a manner that is still quite closely tied to the structure of $R*G$. 

\begin{definition}
\label{defn:equiv} Let $k$ be a field, $R'$ a finitely generated $k$-algebra that is a domain, and let $R$ be a subalgebra of $R'$ generated by $l$-th powers of a set of generators of $R'$ for some $l \in \mathbb{Z}^{>0}$. Suppose that $\varphi_{G}: G \rightarrow \Aut(R')$, for $G$ a finitely generated abelian group, is an action of $G$ on $R'$ that restricts to an action on $R$ with trivial kernel. Let $\Theta$ be a finite set of generators for $G$ closed under taking inverses, and let $h_{G}: \Theta \rightarrow R'$ be a map. Then $A$ is an \textbf{equivariant $k$-algebra} with respect to the data $(R',R,G,\varphi_{G},h_{G})$ if $A$ is isomorphic to the subalgebra of $R'*G$ generated by $R$ and $\{h_{G}(g)e_{g}: g \in \Theta\}$.
\end{definition}

\begin{remark} 
\label{rmk:equiv:faithflat} By construction, $R'$ is a finitely generated $R$-module, hence the corresponding map on prime spectra $\mathbf{p}: \Spec R' \rightarrow \Spec R$ is closed (\cite{Har77}, II Ex. 3.5). Moreover, $\mathbf{p}$ is easily verified to be surjective in this case. A geometric argument for this fact would be as follows: all fibres of $\mathbf{p}$ have the same size (counting multiplicities) and so $\mathbf{p}$ is a flat morphism, and hence open (\cite{Mil80}, Theorem 2.9; \cite{Har77}, III, Ex. 9.1). But $\Spec R'$ is irreducible, hence connected, so the image of $\mathbf{p}$ must be all of $\Spec R$. 
\end{remark} We illustrate Definition \ref{defn:equiv} with a few examples of algebras over a field of characteristic $0$. In particular, we will see that some small quantum groups are examples of equivariant algebras. 

\begin{example}[$\mathcal{O}_{q}((k^{\times})^{2})$]
\label{ex:equiv:qtor} In this example we consider the \textbf{quantum $2$-torus} $\mathcal{O}_{q}((k^{\times})^{2})$, namely the $k$-algebra generated by $U$ and $V$ subject to the relation
\[ UV = qVU \] with $U$ and $V$ also invertible. Here, $k^{\times}$ denotes the multiplicative group of $k$ and $q$ is a non-torsion element in $k^{\times}$. Put $R' = R = k[V^{\pm 1}]$. Then taking $G = \mathbb{Z}$ the integers with the action of $G$ on $R$ defined by ${}^{1}V = qV$, we see that $\mathcal{O}_{q}((k^{\times})^{2})= k\langle R, e_{1} \rangle = R*G$. If $q = \epsilon$ is a primitive $n$-th root of unity in $k$ (and assuming that $k$ contains all $n$-roots of unity) we take $G = \mathbb{Z}/n\mathbb{Z}$ and the same action on $R$ (namely ${}^{1}r = \epsilon r$) also allows us to conclude that $\mathcal{O}_{\epsilon}((k^{\times})^{2}) = R*G$. 
\end{example}

\begin{example}[$U_{q}(\mathfrak{sl}_{2}(k))$]
\label{ex:equiv:qsl2}
Let $q \in k^{\times}$ be a non-torsion element. The \textbf{quantized enveloping algebra} of $sl_{2}(k)$, denoted $U_{q}(sl_{2}(k))$, is defined to be the $k$-algebra with generators $E,F, K^{\pm 1}$ subject to the following relations
\[ KEK^{-1} = q^{2}E \qquad KFK^{-1} = q^{-2}F \qquad EF - FE = \frac{K - K^{-1}}{q - q^{-1}} \] along with $KK^{-1} = K^{-1}K = 1$. 
\\
\\
Take $R' = k[X^{\pm 1}]$, $R = k[X^{\pm 2}]$. Let $G = \mathbb{Z}$ and suppose that $G$ acts on $R'$ by ${}^{1}X = qX$. Then $A$ is isomorphic to a subalgebra of $R'*G$ by taking 
\[ K \mapsto X^{2} \qquad F \mapsto -\frac{X + X^{-1}}{q - q^{-1}}e_{1} \qquad E \mapsto \frac{X + X^{-1}}{q - q^{-1}}e_{-1} \] 

\begin{proof} A routine calculation, but we give the details. Firstly, we note that ${}^{1}X^{2} = q^{2}X^{2}$, and it is clear that the relations $FK = q^{2}KF$ and $EK = q^{-2}KE$ are respected in the image of this map. Now
\begin{eqnarray*}
EF \mapsto -\left(\frac{X + X^{-1}}{q - q^{-1}}e_{-1}\right)\left(\frac{X + X^{-1}}{q - q^{-1}}e_{1}\right) = &  -\frac{1}{(q-q^{-1})^{2}} (X + X^{-1})(q^{-1}X + qX^{-1})e_{0} \\ = &  -\frac{1}{(q - q^{-1})^{2}}(q^{-1}K + q^{-1} + q + qK^{-1})e_{0}
\end{eqnarray*} and similarly 
\begin{eqnarray*}
FE \mapsto -\left(\frac{X + X^{-1}}{q - q^{-1}}e_{1}\right)\left(\frac{X + X^{-1}}{q - q^{-1}}e_{-1}\right) = &  -\frac{1}{(q-q^{-1})^{2}} (X + X^{-1})(qX + q^{-1}X^{-1})e_{0} \\ = &  -\frac{1}{(q - q^{-1})^{2}}(qK + q + q^{-1} + q^{-1}K^{-1})e_{0}
\end{eqnarray*} Taking the difference of these expressions shows that $EF - FE \mapsto \frac{K - K^{-1}}{q - q^{-1}}e_{0}$ as required. 
\end{proof} As in Example \ref{ex:equiv:qtor}, if $q = \epsilon$ is a primitive $n$-th root of unity, then replacing $G$ with $\mathbb{Z}/n\mathbb{Z}$ in the above will give that $U_{\epsilon}(sl_{2}(k))$ is also equivariant. 
\end{example}

\begin{example}
\label{ex:equiv:qSL2} $\mathcal{O}_{q}(SL_{2}(k))$ is defined to be the $k$-algebra generated by $A,B,C,D$ subject to the relations 
\[ AB = qBA \qquad AC = qCA \qquad BC = CB \] \[ BD = qDB \qquad CD = qDC \qquad AD - DA = (q - q^{-1})BC \qquad AD - qBC = 1\] As in previous examples, $q \in k^{\times}$ is non-torsion, and suppose that $q$ has a square root $q^{1/2}$ in $k$. We shall consider a localization of this algebra in which $B$ and $C$ are assumed invertible. Take $R' = k[X^{\pm 1},Y^{\pm 1}]$, $R = [X^{\pm 2}, Y^{\pm 2}]$. Again, we take $G  = \mathbb{Z}$ with its action defined on $R'$ by ${}^{1}X = q^{1/2}X$ and ${}^{1}Y = q^{1/2}Y$. Then a similar calculation to that of Example \ref{ex:equiv:qsl2} will give that $A$ is a subalgebra of $R'*G$ via the homomorphism taking
\[ A \mapsto XYe_{1} \qquad B \mapsto X^{2} \qquad C \mapsto Y^{2} \qquad D \mapsto (XY + qX^{-1}Y^{-1})e_{-1} \] and, as in previous examples, we can also conclude that the analogous localization of $\mathcal{O}_{\epsilon}(SL_{2}(k))$ is an equivariant algebra when $\epsilon$ is a primitive $n$-th root of unity. 
\end{example}

\section{Structures Associated to Equivariant Algebras and their Theories} 
We will now associate a structure to each equivariant $k$-algebra $A$ where $k$ is an algebraically closed field of characteristic $0$. We will carry out the construction, somewhat informally, for a specific example before giving a formal axiomatization of the theory of such a structure in the general case. The associated theory will be shown to be consistent in Proposition \ref{prop:equiv:thycons}. 

\subsection{An example: $U_{q}(\mathfrak{sl}_{2}(k))$}
\label{sect:qsl2}
Recall from Example \ref{ex:equiv:qsl2} that there is an isomorphism of $U_{q}(\mathfrak{sl}_{2}(k))$ with a subalgebra of $R'*G$ given by
\[ K \mapsto X^{2} \qquad F \mapsto -\frac{X + X^{-1}}{q - q^{-1}}e_{1} \qquad E \mapsto \frac{X + X^{-1}}{q - q^{-1}}e_{-1} \] where $G = \mathbb{Z}$, $R' = k[X^{\pm 1}]$ and the action of $G$ on $R'$ is defined by ${}^{1}X = qX$. This isomorphism is used to define a structure to associate to $U_{q}(\mathfrak{sl}_{2}(k))$. We begin by associating a structure that we call a {\em line space} to $U_{q}(\mathfrak{sl}_{2}(k))$. 

\begin{definition}
\label{defn:lsp} A \textbf{line space} is a two-sorted structure $(k,L,\pi, \mathbf{E},C)$ where
\begin{itemize}
\item the sort $k$ is an algebraically closed field of characteristic $0$ endowed with the language of rings;
\item $C$ is a set of constants from $k$; 
\item $\pi: L \rightarrow V(k)$ is a surjective map, where $V$ is some variety defined using constants from $C$ (considered as a definable subset of some cartesian power of $k$); 
\item Fibrewise, $L$ has the structure of a one-dimensional $k$-vector space; namely, the language on $L$ has symbols $\cdot: k \times L \rightarrow L$ and $+: L \times L \rightarrow L$ that restrict to graphs of scalar multiplication and addition on $\pi^{-1}(x)$ for each $x \in V(k)$. 
\item $\mathbf{E} \subseteq L \times V$ is a relation such that for some fixed $l \in \mathbb{Z}^{>0}$,
\begin{itemize}
\item for each $x \in V(k)$, $|\mathbf{E}(L,x)| = l$ and
\item there is a free and transitive action of $\mu_{l}$ on $\mathbf{E}(L,x)$, where $\mu_{l}$ is the group of $l$-th roots of unity in $k$ and the action of $\mu_{l}$ on $\mathbf{E}(L,x)$ is given by scalar multiplication in the fibre $\pi^{-1}(x)$. 
\end{itemize}
\end{itemize}
\end{definition} The parameter $l \in \mathbb{Z}^{>0}$ will be that used to define $R$ as a subalgebra of $R'$; thus for $U_{q}(\mathfrak{sl}_{2}(k))$, $l = 2$. The variety $V$ is that associated to $\Spec R$, hence $V(k) = k^{*} = k \setminus \{0\}$. We will also have reason to utilize the variety $V'$ corresponding to $\Spec R'$; here $V'(k) = V(k)$ and the surjective map $\mathbf{p}: V'(k) \rightarrow V(k)$ is just $\mathbf{p}(x) = x^{2}$. Though $l$ and the variety $V$ form part of the data of a line space, we have not indicated this dependence in the notation $(k,L,\pi,\mathbf{E})$. In what follows, the elements of $\mathbf{E}(L,x)$ for given $x \in k^{*}$ will be referred to as `basis elements' of the fibre $\pi^{-1}(x)$, and we shall drop $\mathbf{E}$ and $C$ from the notation for a line space when discussing them below. 
\\
\\
The following remarks show that there is a representation of $U_{q}(\mathfrak{sl}_{2}(k))$ in an object closely associated with $L$. 

\begin{remark}
\label{rmk:qsl2:charbij} There is a bijective correspondence between points of $k^{*}$ and characters on $R$. If $x \in k^{*}$, its corresponding character on $R$ will be denoted by $\chi_{x}$. 
\end{remark} Remark \ref{rmk:qsl2:charbij} is just the Nullstellensatz and the character $\chi_{x}$ is the homomorphism $R \rightarrow R/\mathfrak{m}_{x} \simeq k$ where $\mathfrak{m}_{x}$ is the maximal ideal of $R$ corresponding to $x$. Evidently, Remark \ref{rmk:qsl2:charbij} also holds for characters on $R'$ and we will use the notation $\chi'_{y}: R' \rightarrow k$ for the character corresponding to $y \in k^{*}$.

\begin{remark}
\label{rmk:qsl2:grpact} The action of $G$ on $R'$ will give a (left) group action on $V'(k)$: if $g \in G$ and $y \in V'(k)$ then $gy$ is the point in $V'(k)$ such that $\chi'_{gy}(r) = \chi'_{y}({}^{g^{-1}}r)$ for all $r \in R$. Specifically, given $y \in V'(k) = k^{*}$ we have $1 \cdot y = q^{-1}y$. By Definition \ref{defn:equiv}, we similarly obtain a left action of $G$ on $V(k)$: for $x \in V(k) = k^{*}$ we have that $1 \cdot x = q^{-2}x$. 
\end{remark} 

\begin{remark}
\label{rmk:qsl2:charlift} There exists a choice of a lift $\chi'_{y:y^{2} = x}$ of $\chi_{x}$ to a character of $R'$ for each $x \in k^{*}$ such that the restrictions of these lifts to $R$ commute with the group action on $k^{*}$.
\end{remark}

\begin{proof} Partition $k^{*}$ into orbits of $G$. Given an orbit, choose a representative $x$ and a square root $y$ of $x$. Then the character $\chi'_{y:y^{2} = x}: R' \rightarrow k$ given by $\chi'_{y:y^{2} = x}(X) = y$ restricts to $\chi_{x}$ and we choose the characters $\chi'_{gy}$ for all remaining points of the orbit. 
\end{proof} 

\begin{remark}
\label{rmk:qsl2:charrep} Let $\{\chi'_{y:y^{2} = x}: x \in k^{*}\}$ be a set of characters of $R'$ that are lifts of the characters on $R$, as in Remark \ref{rmk:qsl2:charlift}. Given $x \in k^{*}$, suppose that we have chosen a basis element $v_{gx} \in L_{gx}:= \pi^{-1}(gx)$ for every $g \in G$. Then $M := \bigoplus_{g \in G} L_{gx}$ is an $R'*G$-module under the action 
\[ (re_{g})v_{g'x} = \chi^{'}_{gg'y:y^{2}= x}(r)v_{gg'x} \] extended linearly to each $L_{gx}$. 
\end{remark}

\begin{proof} It suffices to note that 
\begin{eqnarray*}(re_{g}r'e_{g'})v_{g''x} &  = & \chi'_{g'g''y}(r')\chi'_{gg'g''y}(r)v_{gg'g''x} \\  & = & \chi'_{gg'g''y}({}^{g}r')\chi'_{gg'g''y}(r) v_{gg'g''x} \\ & =& \chi'_{gg'g''y}(r{}^{g}r')v_{gg'g''x} \\ & =& (r{}^{g}r'e_{gg'})v_{g''x} 
\end{eqnarray*}
\end{proof} Given that we can obtain a representation of $U_{q}(\mathfrak{sl}_{2}(k))$ in a module obtained from each orbit of $G$ on $k^{*}$, we try to define the linear maps $K,K^{-1},E$ and $F$ on $L$ in accordance with Remark \ref{rmk:qsl2:charrep}.
\\
\\
Given $x \in k^{*}$, $K$ and $K^{-1}$ must act along $L_{x}$ by multiplying a basis element $v_{x}$ by $\chi_{y:y^{2} = x}'(K^{\pm 1}) = \chi_{x}(K^{\pm 1}) = x^{\pm 1}$ and it is clear that the choice of basis element does not matter. The action of $K^{\pm 1}$ on each $L_{x}$ endows it with the structure of an $R$-module and it is immediate that this action is definable, uniformly for all fibres of $\pi$, in $(k,L,\pi)$ with symbols introduced for $K,K^{-1}$. For $E$ and $F$, Remark \ref{rmk:qsl2:charrep} indicates that we should have 
\[ E: L_{x} \rightarrow L_{q^{2}x} \qquad v_{x} \mapsto \chi'_{qy:(qy)^{2} = q^{2}x}\left(\frac{X + X^{-1}}{q - q^{-1}}\right)v_{q^{2}x}  \] \[ F: L_{x} \rightarrow L_{q^{-2}x} \qquad v_{x} \mapsto \chi'_{q^{-1}y: (q^{-1}y)^{2} = q^{-2}x}\left(-\frac{X + X^{-1}}{q - q^{-1}} \right)v_{q^{-2}x}  \] with respect to a suitable set of characters on $R'$ and a choice of basis elements in each fibre. However, we make the following observation. 

\begin{remark}
\label{rmk:qsl2:nonalgsect} There is no algebraic way of selecting a $G$-equivariant set of characters of $R'$ via a $G$-equivariant section of the map $\mathbf{p}$; namely there is no morphism $s: k^{*} \rightarrow k^{*}$ such that $\mathbf{p} \circ s$ is the identity morphism on $k^{*}$ and $s(gy) = gs(y)$ for all $g \in G$, $y \in k^{*}$. 
\end{remark}

\begin{proof} If such an $s$ did exist, then by Proposition \ref{prop:equiv:thycons} below, the theory associated to $A$ (see Definition \ref{defn:equiv:thy}) would have a model definable in $k$. But this contradicts the result (Corollary \ref{cor:qsl2:nondef}) that no model of this theory is definable in an algebraically closed field. 
\end{proof} Instead, we merely describe $E$ and $F$ in the line space (enriched with symbols for $K,K^{-1},E,F$) by saying that such maps exist with respect to some choice of characters and basis elements, namely by the axiom \[\begin{array}{ll} (\forall v \in L)(\mathbf{E}(v,x) \rightarrow & (\exists v' \in L)(\exists y \in k) \\ & (\mathbf{E}(v',q^{2}x) \wedge y^{2} = x \wedge Ev = \lambda(qy)v' \wedge Fv' = -\lambda(y)v)) \end{array}\] where $\lambda: k^{*} \rightarrow k$ is the function
\[ \lambda(y) = \frac{y+y^{-1}}{q - q^{-1}} \] along with a statement that $E$ and $F$ are $k$-linear on $\pi^{-1}(x)$, and it is the resulting structure that is associated to $U_{q}(\mathfrak{sl}_{2}(k))$. 

\subsection{Associating a structure to an equivariant algebra} 
\label{sect:equiv:thy}
In this subsection, we fix an algebraically closed field $k$ of characteristic $0$ and $A$ an equivariant $k$-algebra with respect to the data $(R',R,G,\varphi_{G},h_{G})$. The theory associated to $A$ will depend on a selection of generators for $R$, the set $\Theta$ and the associated map $h_{G}: \Theta \rightarrow R'$. We enumerate the generators of $A$ corresponding to the map $h_{G}$ as $\{r_{i}e_{g_{i}}: 1 \leq i \leq n\}$ where each $g_{i} \in \Theta$ and $r_{i} = h_{G}(g_{i})$. 
\\
\\
By Remark \ref{rmk:equiv:faithflat}, there is a surjective map $\mathbf{p}: \Spec R' \rightarrow \Spec R$. We can associate a line space $(k,L,\pi)$ to $A$ in the same way as we did for $U_{q}(\mathfrak{sl}_{2}(k))$. Remarks \ref{rmk:qsl2:charbij} to \ref{rmk:qsl2:charrep} (and their proofs) carry over mutatis mutandis to A. The first-order axiomatization of the theory $T_{A}$ of the enriched line space (with additional symbols for the generators of $A$) associated to $A$ is given below. 

\begin{definition}
\label{defn:equiv:thy}
Let 
\[ \mathcal{L}_{A} = (k,L,\pi, \mathbf{E}, C,\mathbf{U}_{i}, H_{j}: 1 \leq j \leq m, 1 \leq i \leq n) \] be the two-sorted language where $(k,L,\pi,\mathbf{E},C)$ is the language of a line space and $\mathbf{U}_{i}$, $H_{j}$ are unary function symbols on $L$ corresponding to the elements $r_{i}e_{g_{i}}$ and generators of $R$ respectively. Then in addition to the axioms of a line space, the theory $T_{A}$ says the following:

\begin{enumerate}
\item \label{defn:equiv:thy:0} For any $1 \leq j \leq m$, 
\[ (\forall x \in V(k))(\forall v \in \pi^{-1}(x))\left(\bigwedge_{j = 1}^{m} H_{j}v = \chi_{x}(H_{j})v \right) \] and $H_{j}$ extends to a linear map on $L_{x} = \pi^{-1}(x)$. 
\item For every $x \in V(k)$ and basis element $v \in L_{x}$, 
\begin{enumerate}
\item \label{defn:equiv:thy:1}for every $1 \leq i \leq n$, there exists $v' \in L$ and $y \in k$ such that
\[ \tilde{\mathbf{U}}_{i}(v,v',y) := \mathbf{E}(v',g_{i}x) \wedge \mathbf{p}(y) = x \wedge \mathbf{U}_{i}v = \chi_{y}'({}^{g_{i}^{-1}}r_{i})v')\] holds, and $\mathbf{U}_{i}$ extends to a linear map on $L_{x}$; 
\item \label{defn:equiv:thy:2} for every $1 \leq i,i' \leq n$,
\[ (\forall v' \in L)(\forall y \in k)(\tilde{\mathbf{U}}_{i}(v,v',y) \rightarrow (\exists v'' \in L)(\exists y' \in k)(\tilde{\mathbf{U}}_{i'}(v',v'',y') \wedge y' = g_{i}y)) \] and;
\item \label{defn:equiv:thy:3} for any $g \in G$ with $g = \prod_{h = 1}^{p} g_{i_{h}}$ as products of elements $g_{i_{h}} \in \Theta$, put 
\[ \mathbf{U}_{\mathbf{h}} := \mathbf{U}_{i_{1}}\dots \mathbf{U}_{i_{p}} \qquad \chi'_{y,\mathbf{h}}(g) := \chi'_{y} (\prod_{h = 1}^{p} {}^{\prod_{q = 0}^{h - 1}g_{i_{q}}} r_{i_{h}}) \qquad g_{i_{0}} := 1_{G} \] and for $\sigma \in S_{p}$ (the permutation group on $p$ letters), put $\mathbf{U}_{\sigma(\mathbf{h})} := \mathbf{U}_{i_{\sigma(1)}}\dots \mathbf{U}_{i_{\sigma(p)}}$ with $\chi'_{y,\sigma(\mathbf{h})}(g)$ defined accordingly. Then we have that for every $\sigma \in S_{p}$,
\[ (\forall v',v'' \in L)(\forall y \in k)(\mathbf{U}_{\mathbf{h}}v = \chi'_{y,\mathbf{h}}(g) v' \wedge \mathbf{U}_{\sigma(\mathbf{h})}v = \chi'_{y,\sigma(\mathbf{h})}(g) v'' \rightarrow v' = v'') \]
\end{enumerate}
\end{enumerate} The set of constants $C$ is extended to include names for all parameters used in \ref{defn:equiv:thy:0} and \ref{defn:equiv:thy:1} - \ref{defn:equiv:thy:3}. 
\end{definition} 

\begin{remark} 
\label{rmk:equiv:thy:com}
\begin{enumerate}
\item \label{rmk:equiv:thy:com:1} If $A$ is an equivariant algebra over an arbitrary field of characteristic $0$, then it is still possible to associate a theory to it. Clearly the theory $T_{A'}$ where $A' := k' \otimes_{R} A$ and $k'$ is an algebraically closed field containing $k$, does not depend on the choice of $k'$. Thus we associate $T_{A'}$ to $A$ for any such $A'$. 
\item That a set of $G$-equivariant characters on $R'$ lifting those on $R$ exists is expressed by the axiom schemes \ref{defn:equiv:thy:1} and \ref{defn:equiv:thy:2}. The axiom scheme \ref{defn:equiv:thy:3} ensures that given an $x \in V(k)$ the basis elements in the fibres $L_{gx}$ for all $g \in G$ are chosen so that $\bigoplus_{g \in G} L_{gx}$ is indeed a representation of $A$ with the actions of $\mathbf{U}_{i}$ defined. For this to be the case, given any basis element $v \in L_{x}$, any two maps taking $v$ to the same target fibre must give scalar multiples (in accordance with Remark \ref{rmk:qsl2:charrep}) of the same target basis element. By the assumption that the action $\varphi_{G}$ restricts to an action with trivial kernel on $R$ in Definition \ref{defn:equiv}, it follows that the action of $G$ on the variety $V(k)$ is free. Thus any fibral coincidences of this kind depend only on the structure of the group $G$ and are independent of the choice of $x \in V(k)$.
\end{enumerate}
\end{remark} The following proposition establishes that a model of $T_{A}$ exists and in the sequel, we will adopt the notation $(k,L)$ for a model of $T_{A}$. 

\begin{proposition}
\label{prop:equiv:thycons} Let $\mathbf{p}: V'(k) \rightarrow V(k)$ be the morphism of varieties corresponding to $\Spec R' \rightarrow \Spec R$. Then there is a model of $T_{A}$ definable in $k$ endowed with a unary function symbol $s$ giving a (set-theoretic) $G$-equivariant section of $\mathbf{p}$. 
\end{proposition}

\begin{proof} It is worth emphasizing that $s$ is not stipulated to be a morphism of varieties. Thus $s: V(k) \rightarrow V'(k)$ is merely a set-theoretic function such that $\mathbf{p} \circ s$ is the identity function on $V(k)$ and $s(gy) = gs(y)$ for all $y \in V'(k)$. We know that such an $s$ exists by the analogue of Remark \ref{rmk:qsl2:charlift} for $A$. 
\\
\\
Let $\mu_{l}$ be the group of $l$-th roots of unity in $k$ and put $\tilde{L} := \mu_{l} \times V(k)$. For each generator $H_{j}$ of $R$ we define the relation $\tilde{H}_{j}$ on $\tilde{L}^{2} \times k$ by 
\[ \tilde{H}_{j}((\gamma,x),(\delta,x'),\alpha) \Leftrightarrow (x = x' \wedge \alpha = \gamma\delta^{-1}\chi_{x}(H_{j}))\] Let $\chi'_{x}$ denote the lift of the character $\chi_{x}$ given by $s$. Then for each element $r_{i}e_{g_{i}}$ of $A$, we introduce a relation $\tilde{\mathbf{U}}_{i}$ on $\tilde{L}^{2} \times k$ where 
\[ \tilde{\mathbf{U}}_{i}((\gamma,x),(\delta,x'),\alpha) \Leftrightarrow (g_{i}x = x' \wedge \alpha = \gamma\delta^{-1}\chi'_{gx}(r_{i}))\] Consider the following equivalence relation on $k \times \tilde{L}$:
\[ (\alpha_{1}, \delta_{1},x_{1}) \sim (\alpha_{2}, \delta_{2}, x_{2}) \Leftrightarrow (\exists \gamma \in \mu_{n})(\alpha_{2} = \gamma\alpha_{1} \wedge  \delta_{2} = \gamma^{-1}\delta_{1}) \] and let $L:= k \times \tilde{L}/\sim$. We shall denote the equivalence class of $(\alpha,\gamma,x)$ in this quotient by $\overline{(\alpha,\gamma,x)}$. Note that there is a natural projection map $\pi: L \rightarrow V(k)$ taking $\overline{(\alpha, \gamma,x)}$ to $x$. 
\\
\\
\textbf{Claim}: Each $L_{x}:= \pi^{-1}(x)$ for $x \in V(k)$ has the structure of a one-dimensional $k$-vector space by
\[ \overline{(\alpha_{1},\delta_{1})} + \overline{(\alpha_{2},\delta_{2})} := \overline{(\gamma^{-1} \alpha_{1} + \alpha_{2}, \delta_{2})} \mbox{ where } \delta_{2} = \gamma \delta_{1} \] \[ \lambda \overline{(\alpha,\delta)} := \overline{(\lambda \alpha,\delta)} \] 

\begin{proof} A routine verification, but we give the details. Suppose that $(\alpha_{1}, \delta_{1}) \sim (\alpha_{1}',\delta_{1}')$ and $(\alpha_{2},\delta_{2}) \sim (\alpha_{2}',\delta_{2}')$ and that $\delta_{2} = \gamma \delta_{1}$. There are $\gamma_{1}, \gamma_{2}$ such that $\delta_{1}' = \gamma_{1} \delta_{1}$ and $\delta_{2}' = \gamma_{2} \delta_{2}$. Thus
\[ \delta_{2}' = \gamma_{2}\gamma\gamma_{1}^{-1} \delta_{1}' \] So it remains to prove that
\[ (\gamma^{-1}\alpha_{1} + \alpha_{2},\delta_{2}) \sim (\gamma_{1}\gamma^{-1}\gamma_{2}^{-1}\alpha_{1}' + \alpha_{2}',\delta_{2}') \] But $\gamma_{1}^{-1}\alpha_{1} = \alpha_{1}'$ and $\gamma_{2}^{-1}\alpha_{2} = \alpha_{2}'$. So
\[ \gamma_{1}\gamma\gamma_{2}^{-1}\alpha_{1}' + \alpha_{2}' = \gamma_{2}^{-1}(\gamma^{-1}\alpha_{1} + \alpha_{2}) \] as required. That scalar multiplication is well-defined is trivial. 
\end{proof} The basis elements in a fibre $L_{x}$ are designated to come from $\mu_{l} \times \{x\} \subseteq \tilde{L}$, hence they are of the form $\overline{(1,\gamma,x)}$ for $\gamma \in \mu_{l}$. By the above claim, there is certainly a free and transitive action of $\mu_{l}$ on the basis elements of $L_{x}$ given by scalar multiplication. The relations $\tilde{H}_{j}$ and $\tilde{\mathbf{U}}_{i}$ allow us to define linear maps on $L$, namely for each $j$
\[ H_{j}: \overline{(1,\gamma,x)} \mapsto \overline{(\alpha,\delta,x')} \Leftrightarrow \tilde{H}_{j}((\gamma,x),(\delta,x'),\alpha) \] is extended $k$-linearly, and similarly for the $\mathbf{U}_{i}$. It is immediate, by construction, that the axioms of Definition \ref{defn:equiv:thy} hold. 
\end{proof} 

\section{Non-algebraicity}
We continue with the notation of Section \ref{sect:equiv:thy}. If $A$ is an equivariant $k$-algebra (with respect to commutative algebras $R,R'$ and the group $G$) and there is a morphism $s: \Spec R \rightarrow \Spec R'$ giving a $G$-equivariant section of $\mathbf{p}: \Spec R' \rightarrow \Spec R$ then Proposition \ref{prop:equiv:thycons} will give that there is a model of $T_{A}$ definable in $k$. In particular those $A$ for which $R = R'$ (for example $A = \mathcal{O}_{q}((k^{\times})^{2})$ of Example \ref{ex:equiv:qtor}) will have this property. 
\\
\\
The question is therefore raised as to what obstructions exist with regard to being able to define or interpret models of $T_{A}$ in an algebraically closed field for a given $A$. Two necessary algebraic conditions for interpretability to be possible are given. The first of these conditions (torsion) is on the map $h_{G}$ forming part of the data of $A$ and is rather stringent. We shall see later that torsion is the algebraic condition characterizing relative categoricity (Theorem \ref{thm:equivthy:cat}), so from a model-theoretic perspective, it is the torsion equivariant algebras that are more interesting. For the torsion equivariant algebras, we obtain a second necessary and sufficient condition for interpretability via a straightforward adaptation of the methods (due to Sustretov) of Sections $3$ and $4$ of the paper \cite{SSZ14}. Using this condition, it can be shown that $U_{q}(\mathfrak{sl}_{2}(k))$ and the localization of $\mathcal{O}_{q}(SL_{2}(k))$ in Example \ref{ex:equiv:qSL2} have associated theories with no model definable in an algebraically closed field. 

\subsection{Interpretability and torsion} 
Firstly, we set up some notation. Let $\mathcal{L},\mathcal{L}'$ be languages, $\mathcal{M}$ an $\mathcal{L}$-structure, $\mathcal{M}'$ a $\mathcal{L}'$-structure. The structure $\mathcal{M}$ is {\em interpretable} in $\mathcal{M}'$ if there is a definable set $M(M')$ in $\mathcal{M'}^{eq}$ corresponding to the universe $M$ of $\mathcal{M}$, and for each predicate $S$ of $\mathcal{L}$, there is an $\mathcal{L}'^{eq}$-definable relation $S(M')$ in $\mathcal{M}'$ such that the structure $(M(M'),S(M'): S \in \mathcal{L})$ is isomorphic as an $\mathcal{L}$-structure to $\mathcal{M}$. If $\mathcal{M}$ is an $\mathcal{L}$-structure and $k$ is an algebraically closed field, if $\mathcal{M}$ is interpretable in $k$ then it is definable in $k$ (namely we can take $M(k)$ and $S(k)$ to be definable sets in $k$) due to the model-theoretic fact that $k$ has elimination of imaginaries. For this reason we shall use the terms definability and interpretability interchangeably. 

\begin{definition} Let $A$ be an equivariant $k$-algebra with respect to the data $(R',R,G,\varphi_{G},h_{G}:\Theta \rightarrow R')$. $A$ is said to be \textbf{torsion} if $h_{G}^{l}(g) \in R$ for every $g \in \Theta$.
\end{definition} All of the examples discussed in Section \ref{sect:equivalgs} are easily seen to be torsion.

\begin{lemma} 
\label{lem:equiv:tors} Suppose that $A$ is an equivariant $k$-algebra, for $k$ an algebraically closed field of characteristic $0$, that is not torsion. Then for any model $(k',L) \models T_{A}$ with $k'$ uncountable containing $k$ and the generic point of $V(k)$, there is an automorphism $\theta$ of $k'$ fixing $k$ that does not extend to an automorphism of $L$. 
\end{lemma}

\begin{proof} Let $(k',L) \models T_{A}$ where $k'$ is uncountable containing $k$ and the generic point $\xi$ of $V(k)$ and suppose for contradiction that every automorphism of $k'$ fixing the constants extends to one of $L$. By axiom \ref{defn:equiv:thy:1} of Definition \ref{defn:equiv:thy}, given a basis element $v_{\xi} \in \pi^{-1}(\xi)$, there is a basis element $v_{g_{i}\xi}$ of $\pi^{-1}(g_{i}\xi)$ such that
\begin{eqnarray} \label{eqn:tors:1} \mathbf{U}_{i}v_{\xi} = \chi'_{y}({}^{g_{i}^{-1}}r_{i})v_{g_{i}\xi} \end{eqnarray} for some $y$ such that $\mathbf{p}(y) = \xi$. Consider the subfield $k(\xi)$ of $k'$ generated by $k$ and $\xi$. If $y'$ is such that $\mathbf{p}(y') = \xi$ with $y \neq y'$, then $\tp^{k'}(y/k(\xi)) = \tp^{k'}(y'/k(\xi))$, where $\tp^{k'}$ denotes the complete type in the language of the field $k'$. By saturation of $k'$, there is an automorphism $\theta: k' \rightarrow k'$ taking $y$ to $y'$.
\\
\\
Let $\tilde{\theta}$ be an extension of $\theta$ to an automorphism $\tilde{\theta}$ of $L$ and let $\{g_{i}: 1 \leq i \leq n\}$ be an enumeration of $\Theta$ with $r_{i} = h_{G}(g_{i})$ for every $i$. Applying $\tilde{\theta}$ to \ref{eqn:tors:1} we obtain that 
\[ \mathbf{U}_{i}\tilde{\theta}(v_{\xi}) = \chi'_{y'}({}^{g_{i}^{-1}}r_{i})\tilde{\theta}(v_{g_{i}\xi}) \]  must hold. The fibres $\mathbf{p}^{-1}(\xi)$ and $\mathbf{p}^{-1}(g_{i}\xi)$ are fixed by $\tilde{\theta}$. Because $\mu_{l}$ acts transitively on the basis elements of each fibre, there are $\gamma,\delta \in \mu_{l}$ such that $\gamma v_{\xi} = \tilde{\theta}(v_{\xi})$ and $\delta v_{g_{i} \xi} = \tilde{\theta}(v_{g_{i} \xi})$. Hence
\begin{eqnarray*} \delta^{-1} \chi'_{y'}({}^{g_{i}^{-1}}r_{i}) v_{g_{i}\xi} &=& \chi'_{y'}({}^{g_{i}^{-1}}r_{i}) \tilde{\theta}(v_{g_{i}\xi}) \\ &=& \mathbf{U}_{i}\tilde{\theta}(v_{\xi}) \\ & = &  \gamma \mathbf{U}_{i} v_{\xi} \\ &= & \gamma \chi'_{y}({}^{g_{i}^{-1}}r_{i}) v_{g_{i}\xi} 
\end{eqnarray*} so we must have that $\delta^{-1} \chi'_{y'}({}^{g_{i}^{-1}}r_{i}) = \gamma \chi'_{y'}({}^{g_{i}^{-1}}r_{i})$. Raising both sides to the $l$-th power,
\[ \chi'^{l}_{y'}({}^{g_{i}^{-1}}r_{i}) = \chi'_{y'}(({}^{g_{i}^{-1}}r_{i})^{l}) = \chi'_{y'}({}^{g_{i}^{-1}}(r_{i}^{l})) = \chi'_{y}({}^{g_{i}^{-1}}(r_{i}^{l})) \] Thus $\chi'_{y}({}^{g_{i}^{-1}}(r_{i}^{l}))$ is Galois-invariant and must therefore lie in $k(\xi) = \Frac(R)$, the field of fractions of $R$. It follows that for some $s_{i}, t_{i}' \in R$, $s_{i}/t_{i}'$ and ${}^{g_{i}^{-1}} r_{i}^{l}$ define the same rational functions on $\mathbf{p}^{-1}(U)$ for $U$ an open subset of $V(k)$. But then $s_{i}$ and ${}^{g_{i}^{-1}} r_{i}^{l}t'_{i}$ give regular functions on $V'(k)$ that agree everywhere (they agree on a dense open set), hence ${}^{g_{i}^{-1}}r_{i}^{l}t'_{i} \in R \Leftrightarrow r_{i}^{l}t_{i} \in R$ where $t_{i} = {}^{g_{i}}t'_{i}$. By the definition of $R'$, if $X$ is a generator in a monomial in $r_{i}^{l}$ with exponent $m$, then every monomial in $t_{i}$ can only contain $X$ raised to an exponent that is a multiple of $l$. But then $m$ must itself be divisible by $l$, so $r_{i}^{l} \in R$ as required. 
\end{proof} 

\begin{lemma} 
\label{lem:fields:defisom} Let $k$ be an infinite field and $K$ an algebraically closed field, both considered as structures in the language of rings. If $k$ is interpretable in $K$ then there is a bijection between $k(K)$ and $K$, definable in $K$, that gives an isomorphism of fields. 
\end{lemma}

\begin{proof} \cite{Pil02}, Theorem 4.13.
\end{proof}

\begin{proposition}
\label{prop:equiv:tors:nondef} If $T_{A}$ has a model definable in an algebraically closed field, then $A$ is torsion. 
\end{proposition}

\begin{proof} Suppose that $(k',L) \models T_{A}$ is interpretable in an algebraically closed field $K$. Then by Lemma \ref{lem:fields:defisom}, $k'(K)$ is definably isomorphic to $K$, so we may assume that $k'(K)$ is interpreted as $K$ with the field operations those given by $K$. Let $K'$ be an uncountable algebraically closed field containing $K$ and the generic point $\xi$ of $V(k') = V(K)$. Then any automorphism fixing $K$ extends to an automorphism of $K'$ and thus induces an automorphism of $L(K')$. So by Lemma \ref{lem:equiv:tors}, $A$ is torsion.
\end{proof}

\subsection{Principal homogeneous spaces} 
The case where $A$ is torsion requires a more detailed analysis; the obstructions to definability in this case are Kummer-theoretic in nature. In this subsection, we introduce the relevant Galois cohomological preliminaries. 

\begin{definition} 
\label{defn:Gprin} Let $\mathcal{G}$ be a group. A group $\mathcal{A}$ is a \textbf{$\mathcal{G}$-group} if there is a left action of $\mathcal{G}$ on $\mathcal{A}$ that is compatible with the group operation on $\mathcal{A}$, i.e. given $\sigma \in \mathcal{G}$ we have ${}^\sigma(ab) = {}^\sigma a{}^\sigma b$ for all $a,b \in \mathcal{A}$. A \textbf{principal homogeneous space} $\mathcal{P}$ for $\mathcal{A}$ is a set endowed with a left action of $\mathcal{G}$ and a right action of the group $\mathcal{A}$ that is free and transitive, and is compatible with the action of $\mathcal{G}$, i.e. for all $\sigma \in \mathcal{G}$, $x \in \mathcal{P}$ and $a \in \mathcal{A}$, ${}^\sigma(x \cdot a) = {}^\sigma x \cdot {}^\sigma a$.
\end{definition} 

\begin{definition}
\label{defn:cohomgrps} Let $\mathcal{G}$ be a group, $\mathcal{A}$ a $\mathcal{G}$-group. We define the group
\[ H^{0}(\mathcal{G},\mathcal{A}) := \{a \in \mathcal{A}: {}^\sigma a = a \mbox{ for all $a \in A$} \} \] and the set
\[ H^{1}(\mathcal{G},\mathcal{A}) := \{h: \mathcal{G} \rightarrow \mathcal{A}: h(\sigma\tau) = h(\sigma){}^\sigma h(\tau) \mbox{ for all $\sigma,\tau \in \mathcal{G}$}\}/\sim \] where $h \sim h'$ if and only if there is $a \in \mathcal{A}$ such that $h(\sigma) = {}^\sigma ah'(\sigma)a^{-1}$. The maps $h$ are called \textbf{cocycles} and two cocycles $h,h'$ such that $h \sim h'$ are said to be \textbf{cohomologous}. 
\end{definition} If $\mathcal{A}$ is abelian then the set $H^{1}(\mathcal{G}, \mathcal{A})$ is naturally a group (with the group operation defined pointwise on elements of $\mathcal{G}$). Definitions \ref{defn:Gprin} and \ref{defn:cohomgrps} both hold in the category of algebraic groups, with the group operations and actions replaced by morphisms in this category. 
\\
\\
Now let $k$ be an arbitrary field, $K/k$ a finite Galois extension, $\mathcal{A}$ an algebraic group defined over $k$. There is a natural action of the Galois group $\Gal(K/k)$ on the $K$-points of $\mathcal{A}$, $\mathcal{A}(K)$, endowing the latter object with the structure of a $\Gal(K/k)$-group. 

\begin{proposition} 
\label{prop:H1prin} 
Let $K/k$ be a finite Galois extension, $\mathcal{A}$ an algebraic group defined over $k$. Then there is a bijective correspondence between the elements of $H^{1}(\Gal(K/k), \mathcal{A}(K))$ and $k$-equivalence classes of principal homogeneous spaces for $\mathcal{A}$ defined over $K$. 
\end{proposition}

\begin{proof} If $\mathcal{P}$ is a principal homogeneous $\mathcal{A}$-space then fixing some $p \in \mathcal{P}(K)$, one obtains a cocycle $h: \Gal(K/k) \rightarrow \mathcal{A}(K)$ given by $h(\sigma) = a_{\sigma}$ where $a_{\sigma}$ is the unique element in $\mathcal{A}(K)$ such that $\sigma p = p \cdot a_{\sigma}$. A different choice of $p$ gives a cohomologous cocycle. The converse is given by Proposition 3.4 of \cite{SSZ14}. 
\end{proof}

If $\mathcal{A}$ is an abelian algebraic group, then there is a group law on the set of principal homogeneous $\mathcal{A}$-spaces and we denote this group by $\WC(K/k, \mathcal{A})$. If $\mathcal{A}$ is a zero-dimensional, $\WC(K/k, \mathcal{A})$ can be described explicitly as follows. 

\begin{definition}
\label{defn:WC} Let $K/k$ be a finite Galois extension, $\mathcal{A}$ be a zero-dimensional algebraic group defined over $k$, $\WC(K/k,\mathcal{A})$ the set of principal homogeneous $\mathcal{A}$-spaces defined over $K$. For $\mathcal{P},\mathcal{P}' \in \WC(K/k,\mathcal{A})$, we define $\mathcal{P}*\mathcal{P}'$ to be the quotient of the direct product $\mathcal{P} \times \mathcal{P}'$ by the action $(x,y) \cdot a = (x\cdot a,y \cdot a^{-1})$ for $a \in \mathcal{A}(K)$: the quotient is a principal homogeneous space via the action $(x,y)\cdot a = (x\cdot a,y)$. The inverse $\mathcal{P}^{-1}$ is defined to be the opposite of $\mathcal{P}$. These two operations endow $\WC(K/k,\mathcal{A})$ with the structure of a group. 
\end{definition}

\begin{remark}
\label{rem:WC:H1bij} It can be verified that the bijection in Proposition \ref{prop:H1prin} extends to a group isomorphism $\WC(K/k,\mathcal{A}) \simeq H^{1}(\Gal(K/k), \mathcal{A}(K))$. Under this isomorphism, if $h,h'$ are cocycles with corresponding principal homogeneous spaces $\mathcal{P}_{h}, \mathcal{P}_{h'}$, then $\mathcal{P}_{h}*\mathcal{P}_{h'} = \mathcal{P}_{hh'}$ and $\mathcal{P}^{-1}_{h} = \mathcal{P}_{h^{-1}}$. 
\end{remark} We remark that all of the above extends to infinite Galois extensions. In particular $H^{1}(\Gal(k^{sep}/k),\mathcal{A}(k^{sep}))$ classifies principal homogeneous $\mathcal{A}$-spaces defined over any Galois extension of $k$ (see \cite{Ser64}, \S 5.2).

\subsection{A criterion for non-definability}
We now give the analogue of \cite{SSZ14}, Theorem 4.4 in the current setting after stating the main theorem of Kummer theory. 

\begin{theorem}
\label{thm:kum} Let $k$ be a perfect field, $l$ an integer that does not divide the characteristic of $k$. Suppose that $k$ contains the group of $l$-th roots of unity. Then 
\[ H^{1}(\Gal(k^{sep}/k),\mu_{l}) \simeq k^{\times}/(k^{\times})^{l} \] 
\end{theorem}

\begin{proof} \cite{Ser64}, II \S 1.2. 
\end{proof}

\begin{remark}
\label{rem:kum} The isomorphism $\WC(k^{sep}/k, \mu_{l}) \simeq k^{\times}/(k^{\times})^{l}$ given by (the infinite Galois-theoretic analogue of) Proposition \ref{prop:H1prin} and Theorem \ref{thm:kum} can be described explicitly. If $x \in k^{\times}$ then there is a corresponding principal homogeneous space for $\mu_{l}$ given by $y\mu_{l}$ for some $l$-th root $y$ of $x$. Conversely, because the action of $\Gal(k^{sep}/k)$ on $\mu_{l}$ is trivial, a principal homogeneous space corresponds to a homomorphism $h: \Gal(k^{sep}/k) \rightarrow \mu_{l}$. Then $\mu_{l}/\ker h$ is the Galois group of a cyclic field extension $k'/k$ of exponent dividing $l$. If $\sigma$ generates this Galois group and $h(\sigma) = \epsilon$ (for $\epsilon$ a primitive $l$-root of unity), then there is $y \in k'$ such that $\sigma(y) = \epsilon y$. It follows that $y^{l}$ is invariant under the galois group, so $y^{l} \in k$. 
\end{remark}

\begin{theorem}
\label{thm:equiv:non-def} Let $A$ be a torsion equivariant $k$-algebra for $k$ an algebraically closed field of characteristic $0$ with respect to the data $(R',R,G,\varphi_{G},h_{G})$. Then the theory $T_{A}$ has a model definable in an algebraically closed field if and only if $H_{G}(g) \in R$ for every $g \in \Theta$. 
\end{theorem}

\begin{proof} Firstly we note that because $A$ is torsion, $H_{G}(g)^{l} \in R$ (and is nonzero) for every $g \in \Theta$, so it makes sense to talk about the class of such an element in $\Frac(R)^{\times}/(\Frac(R)^{\times})^{l}$. 
\\
\\
Now suppose that $(k',L) \models T_{A}$ is interpretable in an algebraically closed field $K$ over the subfield $\mathbb{F}$ of $K$. Then (as in the proof of Proposition \ref{prop:equiv:tors:nondef}) we may assume that $k'(K) = K$. Let $K'$ be an algebraically closed field containing $K$ and the generic point $\xi$ of $V(K)$. For every $x \in V(K')$, put $\tilde{L}_{x} := \mathbf{E}(L(K'),x)$. By the axioms for $T_{A}$, $\tilde{L}_{\xi}$ is a principal homogeneous space for $\mu_{l}$, definable over the field $\mathbb{F}(\xi)$. This principal homogeneous space corresponds to a cocycle class in $H^{1}(\Gal(\mathbb{F}(\xi)^{sep}/\mathbb{F}(\xi)),\mu_{l}) \simeq \mathbb{F}(\xi)^{\times}/(\mathbb{F}(\xi)^{\times})^{l}$ by Proposition \ref{prop:H1prin} and Theorem \ref{thm:kum}. Let $g_{i} \in \Theta$, $h_{G}(g_{i}) = r_{i}$. 
\\
\\
\textbf{Claim}: $\tilde{L}_{\xi}*\tilde{L}^{-1}_{g \xi} = \chi'_{\nu}({}^{g_{i}^{-1}}r_{i})\mu_{l}$ for some $\nu$ such that $\mathbf{p}(\nu) = \xi$.

\begin{proof} The relation $\tilde{\mathbf{U}}_{i}$ of Definition \ref{defn:equiv:thy:1} is definable by assumption and hence its restriction to  $\tilde{L}_{\xi} \times \tilde{L}_{g_{i}\xi} \times \{\nu\}$ is also definable for some $\nu$ with $\mathbf{p}(\nu)= \xi$. But for every $\gamma \in \mu_{l}$, if $\tilde{\mathbf{U}}_{i}(v,v',\nu)$ holds for $(v,v') \in \tilde{L}_{\xi} \times \tilde{L}_{g_{i}\xi}$ then so does $\tilde{\mathbf{U}}_{i}(\gamma v, \gamma v', \nu)$. Hence the claim follows by Definition \ref{defn:WC}. 
\end{proof} Both $\tilde{L}_{\xi}$ and $\tilde{L}_{g_{i}\xi}$ must have the same class in $\mathbb{F}(\xi)^{\times}/(\mathbb{F}(\xi)^{\times})^{l}$ because the action of $G$ on $V'(K')$ gives an algebraic isomorphism between them, definable over $\mathbb{F}(\xi)$. So the class of $\chi'_{\nu}({}^{g_{i}^{-1}}r_{i})\mu_{l}$, which is ${}^{g_{i}^{-1}}r_{i}^{l} \mbox{  $\mod$  }  (\Frac(R)^{\times})^{l}$, is trivial. Hence there exist $s_{i},t_{i} \in R^{\times}$ such that ${}^{g_{i}^{-1}}r_{i}^{l}s_{i}^{l} = t_{i}^{l}$, which holds if and only if ${}^{g_{i}^{-1}}r_{i}s_{i} = t_{i}$ (in $R'$). Now we conclude (as in the proof of Lemma \ref{lem:equiv:tors}) that the generators of any monomial of ${}^{g_{i}^{-1}}r_{i}$ must have exponent divisible by $l$, so ${}^{g_{i}^{-1}}r_{i}$ (and hence $r_{i}$) is in $R$. The converse is immediate. 
\end{proof}

\begin{corollary}
\label{cor:qsl2:nondef} No model of $T_{A}$ for $A = U_{q}(\mathfrak{sl}_{2}(k))$ and the localization of $\mathcal{O}_{q}(SL_{2}(k))$ considered in Example \ref{ex:equiv:qSL2} is definable in an algebraically closed field. 
\end{corollary}

\section{Categoricity, Quantifier Elimination and Zariski Structure} 
Much of the remainder of this paper consists of a straightforward adaptation of the methods and results of Sections $5$ and $6$ of \cite{SSZ14}. Consequently, a number of proofs are omitted, with references given as appropriate. 

\subsection{Categoricity} 
The theorem in this subsection demonstrates that torsion algebraically characterizes those $A$ for which models of $T_{A}$ are relatively categorical. Before stating the theorem, we give the following definition. 

\begin{definition}
\label{defn:equiv:pathlength} Let $A$ be an equivariant $k$-algebra, $(k',L) \models T_{A}$ and let $\Lambda$ be a set of representatives for the partition of $V(k')$ into orbits of $G$. Given $z \in V(k')$, the \textbf{path-length} of $z$ with respect to $\Lambda$, denoted $l(z)$, is defined to be the minimal word length of $g$ in the generators $\Theta$ such that $z = {}^{g}x$. 
\end{definition}

\begin{theorem} 
\label{thm:equivthy:cat}
Let $A$ be an equivariant $k$-algebra. Then $A$ is torsion if and only if for all $(k',L) \models T_{A}$ with $k'$ uncountable containing $k$ and the generic point of $V(k)$, any automorphism $\theta$ of $k'$ fixing the constants extends to an automorphism $\tilde{\theta}$ of $L$. 
\end{theorem} 

\begin{proof} The implication from right to left is given by Lemma \ref{lem:equiv:tors}. For the converse, let
\[ V(k') = \bigcup_{x \in \Lambda} Gx \] be a partition of $V$, $\Lambda$ a set of representatives. Suppose that $\theta$ is an automorphism of $k'$ fixing the constants. We extend $\theta$ to an automorphism of $L$ by inducting on the path-length of $z \in V(k')$. 

\begin{enumerate}
\item $l(z) = 0$. Then $z = x$ for some representative $x$. Pick any basis elements $v_{x} \in \pi^{-1}(x)$, $v_{\theta x} \in \pi^{-1}(\theta x)$ and define $\tilde{\theta}(v_{x}):= v_{\theta x}$. We then extend $\tilde{\theta}$ linearly to the fibre $\pi^{-1}(x)$. 
\item \label{thm:equivthy:cat:ind1} $l(z) = 1$. Then $z = g_{i}x$ for some $g_{i} \in \Theta$ and representative $x$. By axiom \ref{defn:equiv:thy:1} of Definition \ref{defn:equiv:thy}, given some basis element $v_{x} \in \pi^{-1}(x)$, there is $y \in k$ such that $\mathbf{p}(y) = x$ and 
\begin{eqnarray} (k',L) \models \mathbf{U}_{i} v_{x} = \chi'_{y}({}^{g_{i}^{-1}}r_{i})v_{z} \end{eqnarray} for some basis element $v_{z}$ of $\pi^{-1}(z)$. Correspondingly, for $v'_{x}:= \tilde{\theta}(v_{x}) \in \pi^{-1}(\theta x)$, there is $y' \in k$ such that $\mathbf{p}(y') = \theta x$ and 
\begin{eqnarray} \label{eqn:cat:1} (k',L) \models \mathbf{U}_{i} v'_{\theta x} = \chi'_{y'}({}^{g_{i}^{-1}}r_{i})v'_{\theta z} \end{eqnarray} for some $v'_{\theta z} \in \pi^{-1}(\theta z)$. Now $\mathbf{p}(\theta y) = \theta x$. Because $A$ is torsion, $r_{i}^{l} \in R$. Thus 
\[ \chi'^{l}_{y'}({}^{g_{i}^{-1}}r_{i}) = \chi_{\theta x}({}^{g_{i}^{-1}}(r_{i}^{l})) = \chi'^{l}_{\theta y}({}^{g_{i}^{-1}}r_{i}) \] giving that $\chi'_{y'}({}^{g_{i}^{-1}}r_{i}) = \gamma \chi'_{\theta y}({}^{g_{i}^{-1}}r_{i})$ for some $\gamma \in \mu_{l}$. Hence \ref{eqn:cat:1} gives
\[ (k',L) \models \mathbf{U}_{i} v'_{\theta x} = \gamma \chi'_{\theta y}({}^{g_{i}^{-1}}r_{i})v'_{\theta z} \] and we put $\tilde{\theta}(v_{z}) := \gamma v'_{\theta z}$, extended linearly to $\pi^{-1}(z)$. 
\item $l(z) > 1$. By induction, $\tilde{\theta}$ has already been extended to the fibre $\pi^{-1}(g_{i'}^{-1}z)$ for some $g_{i'} \in \Theta$. Choose any $g_{i} \in \Theta$ and put $w = g_{i}^{-1}g_{i'}^{-1}z$. For some $y \in \pi^{-1}(w)$ we have that 
\begin{eqnarray} \label{eqn:cat:2} (k',L) \models \mathbf{U}_{i}v_{w} = \chi'_{y}({}^{g_{i}^{-1}}r_{i})v_{g_{i}w} \end{eqnarray}  for some basis elements $v_{w} \in \pi^{-1}(w)$ and $v_{g_{i}w} \in \pi^{-1}(g_{i}w)$; moreover the corresponding transform of this equation under $\tilde{\theta}$ also holds. By axiom \ref{defn:equiv:thy:2} of Definition \ref{defn:equiv:thy}, there is a basis element of $v_{z} \in \pi^{-1}(z)$ such that $y' = g_{i}y$ and 
\[ (k',L) \models \mathbf{U}_{i}v_{g_{i}w} = \chi'_{y'}({}^{g_{i'}^{-1}}r_{i'})v_{z} \] But $\theta y' = g_{i} \theta y$, so we extend $\tilde{\theta}(v_{z}):= v_{\theta z}$ where $v_{\theta z}$ is the basis element of $\pi^{-1}(\theta z)$ given by applying axiom \ref{defn:equiv:thy:2} of Definition \ref{defn:equiv:thy} to the $\tilde{\theta}$-transform of \ref{eqn:cat:2}. 
\end{enumerate}
\end{proof}

\subsection{Quantifier elimination} 
For the rest of this paper, we assume that $A$ is torsion. Fix some model $(k',L) \models T_{A}$. We provide some motivation for the definable sets we wish to consider as giving an elimination set. Let $v = (v_{1}, \dots, v_{s})$ be a tuple from the sort $L$. The $v_{i}$ can be re-indexed according to the fibers of $\pi$ in which they appear, i.e. there is an enumeration $(v_{ij}: 1 \leq i \leq t; 1 \leq j \leq s_{i}, \sum_{i} s_{i} = s)$ such that given $v_{ij}, v_{kl}$, we have $i = k$ if and only if $\pi(v_{ij}) = \pi(v_{kl})$. Let $x_{i} \in V(k')$ be such that $\pi(v_{ij}) = x_{i}$ for all $i$. Then given a basis element $v_{i}$ in each fibre $\pi^{-1}(x_{i})$, there are scalars $\lambda_{ij} \in k$ such that 
\[ \models \bigwedge_{i=1}^{t} \bigwedge_{j=1}^{s_{i}} \lambda_{ij} v_{i} = v_{ij} \] One expects that the sentences satisfied by $v$ can be determined from the relationships between the $v_{i}$, and that these relationships are precisely those arising from the maps that exist between the corresponding fibres. Because there are only existential statements in $T_{A}$ governing how the maps moving between different fibres behave, we can only expect to eliminate quantifiers to the level of existential sentences in general. 
\\
\\
We now set up some notation. Suppose that $v_{i}$ and $v_{j}$ for $i < j$ are basis elements lying in fibres over the same orbit of $G$ and that $g_{ij}x_{i} = x_{j}$ for some $g_{ij} \in G$. Then with respect to a fixed presentation of $g$ as a product of elements of $\Theta$, there is an existential sentence coding the corresponding path from $\pi^{-1}(x_{i})$ to $\pi^{-1}(x_{j})$; namely that there exist $y_{ij} \in k$ and $\gamma_{ij} \in \mu_{l}$ such that
\[ \phi^{g_{ij}}_{ij}(v_{i},v_{j},y_{ij},\gamma_{ij}) := (\exists v_{ij} \in L)(\mathbf{U}_{\mathbf{h}}v_{i} = \chi'_{y_{ij},\mathbf{h}}(g)v_{ij} \wedge v_{j} = \gamma_{ij} v_{ij}) \] holds, where $\mathbf{U}_{\mathbf{h}}$ and $\chi'_{y_{ij},\mathbf{h}}(g)$ are as in Definition \ref{defn:equiv:thy} \ref{defn:equiv:thy:3}.

\begin{definition} 
\label{defn:forms:core}
Let $(v_{ij}: 1 \leq i \leq t; 1 \leq j \leq s_{i}, \sum_{i} s_{i} = s)$ and $x = (x_{1}, \dots x_{r})$ be tuples of variables from the sorts $L$ and $k'$ respectively. A \textbf{core formula} with variables $(v,x)$ is defined to be a formula of the following shape:
\[ \exists_{i = 1}^{t} v_{i} \exists \lambda \exists y \exists \gamma \left(\bigwedge_{i=1}^{t} \bigwedge_{j=1}^{s_{i}} \lambda_{ij} v_{i} = v_{ij} \wedge \bigwedge_{(i,j) \in \Xi} \phi_{ij}^{g_{ij}}(v_{i},v_{j},y_{ij},\gamma_{ij}) \wedge S(\lambda,y,\gamma,\pi(v),x)\right)\] where

\begin{enumerate}
\item \label{defn:forms:core:1} $\Xi$ is a subset of $\{(i,j): 1 \leq i < j \leq t\}$; 
\item \label{defn:forms:core:2} $\lambda = (\lambda_{ij}: 1 \leq i \leq t, 1 \leq j \leq s_{i})$, $y = (y_{ij}: (i,j) \in \Xi)$, $\gamma = (\gamma_{ij}: (i,j) \in \Xi)$ and $\pi(v) = (\pi(v_{i}): 1 \leq i \leq t)$;
\item $S$ defines a Zariski constructible subset of $k'^{r_{1}} \times V(k')^{t} \times \mu_{l}^{r_{2}}$ where
\begin{enumerate}
\item $r_{1} = l(x) + l(y) + s + t$ (where $l$ denotes length of the tuple);
\item $r_{2} = l(\gamma)$.
\end{enumerate}
\end{enumerate} A \textbf{core type} is defined to be a consistent collection of core formulas. If $(v,a)$ is a tuple of elements from $L^{s} \times k'^{r}$, $D$ a set of parameters from $k'$, the \textbf{core type of $(v,a)$ over $D$} (denoted $\ctp (v,a/D)$) is defined to be the set of all core formulas satisfied by $(v,a)$ with parameters from $D \cup C$. 
\end{definition} 

\begin{proposition} 
\label{prop:core:quantelim}
Let $(k',L) \models T_{A}$ be $\aleph_{0}$-saturated. Suppose that $(v,c), (w,d)$ are both tuples from $L^{s} \times k'^{r}$ with the property that $\ctp(v,c/D) = \ctp(w,d/D)$. Then $\tp(v,c/D) = \tp(w,d/D)$. 
\end{proposition}

\begin{proof} The proof is similar to that of Theorem \ref{thm:equivthy:cat} insofar as we construct an automorphism $\tilde{\theta}$ of $(k',L)$ that takes $(v,c)$ to $(w,d)$. Re-index the tuple $v$ as $(v_{ij}: 1 \leq i \leq t; 1 \leq j \leq s_{i}, \sum_{i} s_{i} = s)$ so that given $v_{ij}, v_{kl}$, we have $i = k$ if and only if $\pi(v_{ij}) = \pi(v_{kl})$. Then there exists
\begin{itemize}
\item a tuple $x^{(1)} = (x^{(1)}_{i})$ of elements of $k'$ such that $\pi(v_{ij}) = x^{(1)}_{i}$ for every $i$; and 
\item a tuple of basis elements $(v^{(1)}_{i})$ and a tuple of scalars $\lambda^{(1)} = (\lambda^{(1)}_{ij})$ such that 
\[ \models \bigwedge_{i = 1}^{t} \bigwedge_{j=1}^{s_{i}} \lambda_{ij}^{(1)} v_{i}^{(1)} = v_{ij} \] 
\end{itemize} Let $\Xi$ consist of those elements $(i,j)$ for $1 \leq i < j \leq t$ such that there exists $g_{ij} \in G$ with $g_{ij}x^{(1)}_{i} = x^{(2)}_{j}$. Fix a presentation for each $g_{ij}$. Then the formula 
\[ \varphi(v, v^{(1)},\lambda^{(1)}, y^{(1)}, \gamma^{(1)}) := \bigwedge_{i=1}^{t} \bigwedge_{j=1}^{s_{i}} \lambda_{ij} v_{i} = v_{ij} \wedge \bigwedge_{(i,j) \in \Xi} \phi_{ij}^{g_{ij}}(v_{i},v_{j},y_{ij},\gamma_{ij}) \] holds for some tuples $y^{(1)}$ and $\gamma^{(1)}$. Now consider the following set of formulas,
\[ \Sigma := \begin{array}{ll} \{\varphi(w, w',\lambda', y', \gamma') \wedge S(\lambda',y',\gamma',\pi(w'),d): \\ \qquad (k',L) \models \varphi(v, v^{(1)},\lambda^{(1)}, y^{(1)}, \gamma^{(1)}) \wedge S(\lambda^{(1)}, y^{(1)},\gamma^{(1)},x^{(1)}, c) \}  \end{array} \] where the variables have been primed to distinguish them from parameters and the $S$ range over all constructible subsets of a cartesian power of $k$ with parameters from $D$. 
\\
\\
\textbf{Claim}: $\Sigma$ is consistent. 

\begin{proof} It suffices to show that $\Sigma$ is finitely consistent. It is clear that $\Sigma$ is closed under finite conjunctions. If 
\[ (k',L) \models \varphi(v, v^{(1)},\lambda^{(1)}, y^{(1)}, \gamma^{(1)}) \wedge S(\lambda^{(1)}, y^{(1)},\gamma^{(1)},x^{(1)}, c) \] for some $\varphi \in S \in \Sigma$ then quantifying out the $v^{(1)},\lambda^{(1)}, y^{(1)}, \gamma^{(1)}$ gives a core formula over $D$ that is also satisfied by $(w,d)$. So there exist $w^{(2)},\lambda^{(2)}, y^{(2)}, \gamma^{(2)}$ such that 
\[ (k',L) \models \varphi(w, w^{(2)},\lambda^{(2)}, y^{(2)}, \gamma^{(2)}) \wedge S(\lambda^{(2)},y^{(2)},\gamma^{(2)},\pi(w^{(2)}),d) \] as required. 
\end{proof} By saturation of $(k',L)$, there is a tuple $(w^{(2)},\lambda^{(2)}, y^{(2)}, \gamma^{(2)})$ satisfying $\Sigma$. In particular, $\tp^{k'}(\lambda^{(1)}, y^{(1)},\gamma^{(1)},x^{(1)}, c) = \tp^{k'}(\lambda^{(2)},y^{(2)},\gamma^{(2)},\pi(w^{(2)}),d)$ and there is an automorphism $\theta$ of $k'$ taking $(\lambda^{(1)}, y^{(1)},\gamma^{(1)},x^{(1)}, c)$ to $(\lambda^{(2)},y^{(2)},\gamma^{(2)},\pi(w^{(2)}),d)$. It remains to extend $\theta$ to an automorphism $\tilde{\theta}$ of $L$ and to do this we proceed as in the proof of Theorem \ref{thm:equivthy:cat}, with a minor adjustment. Let 
\[ V(k') = \bigcup_{x \in \Lambda} Gx \] be a partition, $\Lambda$ a set of representatives. For $x \in \Lambda$, suppose that the fibres over $Gx$ contain $v^{(1)}_{q_{1}},\dots, v^{(1)}_{q_{m}}$ with the indices ordered so that $q_{l} < q_{l+1}$ for $1 \leq l \leq m$. Then without loss, we can take $x = v^{(2)}_{q_{1}}$. If $g_{q_{1}q_{2}} = \prod_{h=1}^{p} g_{i_{h}}$ for $g_{i_{h}} \in \Theta$ then for $z = g_{i_{1}}x$ we modify the induction step \ref{thm:equivthy:cat:ind1} in the proof of Theorem \ref{thm:equivthy:cat} so that we select $y^{(1)}_{ij}$ instead of an arbitrary lift of $x$ to $V'(k')$ given by axiom \ref{defn:equiv:thy:1} of Definition \ref{defn:equiv:thy}. It is then immediate that $\tilde{\theta}(v^{(1)}_{q_{l}}) = w^{(2)}_{q_{l}}$ for every $l$. 
\end{proof} By compactness, every $\mathcal{L}_{A}$-formula with parameters from the field sort is then equivalent to a boolean combination of core formulas. Some further analysis reveals the structure of subsets of $(k',L)$ defined using parameters from both $L$ and $k'$. 

\begin{definition} 
\label{defn:forms:gencore}
Let $v'$ be a tuple of elements from $L$ with length $p$ such that all elements of $v'$ are basis elements. Let $v = (v_{1}, \dots, v_{m})$, $w = (w_{1}, \dots, w_{n})$ be tuples of variables from $L$. A \textbf{general core formula} with variables $(v,w,x)$ over $v'$ with parameters $D$ from $k'$ is a formula of the following shape:
\[ \exists_{i = 1}^{t} v_{i} \exists \lambda \exists \mu \exists y \exists \gamma \left(\bigwedge_{i=1}^{t} \bigwedge_{j=1}^{s_{i}} \lambda_{ij} v_{i} = v_{ij} \wedge \phi \wedge \bigwedge_{(i,j) \in \Xi} \phi_{ij}^{g_{ij}} \wedge S(\lambda,\mu,y,\gamma,\pi(v),x) \right) \] where 
\begin{enumerate}
\item $(v_{ij}: 1 \leq i \leq t, 1 \leq j \leq s_{i}, \sum s_{i} = s)$ is an enumeration of variables for $v$;
\item $(w_{ij}: 1 \leq i \leq q, 1 \leq j \leq p_{i}, \sum p_{i} = q)$ is an enumeration of variables for $w$;
\item $\phi$ is defined to be 
\[ \bigwedge_{i=1}^{q} \bigwedge_{j=1}^{p_{i}} \mu_{ij} v'_{i} = w_{ij} \wedge \bigwedge_{(i,j) \in \Xi_{1}} \phi_{ij}^{g_{ij}}(v'_{i},v_{j},y_{ij},\gamma_{ij}) \wedge \bigwedge_{(i,j) \in \Xi_{2}} \phi_{ij}^{g_{ij}}(v_{i}, v'_{j}, y_{ij},\gamma_{ij})\] 
where 
\[ \Xi_{1} \subseteq \{(i,j): 1 \leq i \leq q, 1 \leq j \leq t\} \qquad \Xi_{2} \subseteq \{(i,j): 1 \leq i \leq t, 1 \leq j \leq q \} \] 
\item $\Xi$ and $\lambda$ are as in Definition \ref{defn:forms:core} \ref{defn:forms:core:1} and \ref{defn:forms:core:2}, $y = (y_{ij})$, $\gamma = (\gamma_{ij})$ taking indices in $\Xi \sqcup \Xi_{1} \sqcup \Xi_{2}$; 
\item $S$ defines a Zariski constructible subset of $k'^{r_{1}} \times V(k')^{t} \times \mu_{l}^{r_{2}}$ where
\begin{enumerate}
\item $r_{1} = l(x) + l(y) + s + t + q$;
\item $r_{2} = l(\gamma)$.
\end{enumerate}
\end{enumerate}
We shall denote such a formula by $\exists S$ and call $S$ the \textbf{Zariski constructible component} of $\exists S$. 
\end{definition}

\begin{proposition} 
\label{prop:gencore:quantelim} 
Let $(k',L) \models T_{A}$. Any formula with parameters from $L$ and $k'$ is then equivalent to a boolean combination of general core formulas. 
\end{proposition}

\begin{proof} Analogous to the proof of Proposition 5.2 of \cite{SSZ14}. 
\end{proof}

\subsection{Constructibility and Zariski Structure} 
Proposition \ref{prop:gencore:quantelim} suggests taking sets of the form $\exists C$ (where $C$ defines a closed subset of a cartesian power of $k$) as giving the closed sets of a topology on $k'^{m} \times L^{n}$ for some $m,n$. However, there is a priori no guarantee that a given definable set will be constructible for sets definable by formulas of this kind. It transpires that if the $C$ are taken to have a particular form, we do indeed have constructibility. 

\begin{definition}
\label{defn:muinv}
Let $C$ be a formula in the language of the field sort defining a closed subset of $k^{r_{1}} \times V(k')^{t} \times \mu_{l}^{r_{2}}$. We define the \textbf{action} of $\delta \in \mu_{l}^{r_{2}}$ on $C$ to be 
\[ C^{\delta} = \{(\lambda_{ij}, \mu, y, \gamma, z, x): (\delta_{i}^{-1}\lambda_{ij}, \mu, y, \delta \cdot \gamma, z, x) \in C\} \] where 
\[ \delta \cdot \gamma = \left\{ \begin{array}{ll} \delta_{i}^{-1} \gamma_{ij} \delta_{j} & (i,j) \in \Xi \\ \gamma_{ij} \delta_{j} & (i,j) \in \Xi_{1} \\ \delta_{i}^{-1}\gamma_{ij} & (i,j) \in \Xi_{2} \end{array} \right. \] $C$ is defined to be \textbf{$\mu_{l}$-invariant} if $C^{\delta} = C$ for every $\delta \in \mu_{l}^{r_{2}}$. 
\end{definition} 

\begin{proposition} 
\label{prop:const}
All definable subsets of $(k',L)$ are constructible: every definable subset of $(k',L)$ is a boolean combination of those defined by general core formulas $\exists C$ where $C$ is Zariski closed and $\mu_{l}$-invariant. 
\end{proposition} 

\begin{proof} Analogous to \cite{SSZ14}, Proposition 5.3. 
\end{proof} We introduce a topology on $k'^{m} \times L^{n}$ by taking as a basis of closed sets those subsets of $k'^{m} \times L^{n}$ that are defined by general core formulas $\exists C(v,w,x)$ (for $(v,w)$ a tuple of variables from $L^{n}$ and $x \in k^{m}$) where $C$ is Zariski closed and $\mu_{l}$-invariant. Closed sets are given by finite unions and arbitrary intersections of basic closed sets. If $n = 0$, then these formulas reduce to those of the form $C(x)$ where $C$ defines a Zariski closed subset of $k^{m}$, hence the topology on $(k',L)$ gives us the classical Zariski topology on the sort $k'$ and its cartesian powers.

\begin{proposition} 
\label{prop:top:Noeth}
The topology defined on $(k',L)$ is Noetherian. 
\end{proposition}

\begin{proof} See \cite{SSZ14}, Proposition 6.1. 
\end{proof} 

\begin{definition}
\label{defn:dim} Let $\exists C$ define a basic closed irreducible subset of $k'^{m} \times L^{n}$. The \textbf{dimension} of $\exists C(k',L)$ is defined to be the dimension of $C(k')$. For $\exists C$ defining a closed set,
\[ \dim \exists C(k',L):=\max\{C_{i}\} \] where the $C_{i}$ are the irreducible components of $C$. If $\exists S$ is constructible, its dimension is defined to be the dimension of its closure. 
\end{definition}

\begin{theorem} 
\label{thm:Zarstruct}
$(k',L)$ is a Zariski structure which is presmooth if $V(k')$ is smooth.
\end{theorem}

\begin{proof} Analogous to \cite{SSZ14}, Theorem 6.12. 
\end{proof}

\end{document}